\documentclass[english]{article}
\usepackage[T1]{fontenc}
\usepackage[latin9]{inputenc}
\usepackage{array}
\usepackage{amsmath}
\usepackage{graphicx}
\usepackage{amssymb}

\makeatletter

\newcommand{\noun}[1]{\textsc{#1}}
\providecommand{\tabularnewline}{\\}

\newcommand{\lyxaddress}[1]{
\par {\raggedright #1
\vspace{1.4em}
\noindent\par}
}

\providecommand{\tabularnewline}{\\}

\usepackage{babel}

\makeatother

\usepackage{babel}

\makeatother

\usepackage{babel}

\makeatother

\usepackage{babel}

\makeatother

\usepackage{babel}

\begin{document}

\title{An Embedded Boundary Method for Two Phase Incompressible Flow}

\author{S. Wang, J. Glimm, R. Samulyak, X. Jiao, and C. Diao}

\maketitle

\lyxaddress{Department of Applied Mathematics and Statistics, Stony Brook University,
Stony Brook, NY 11794 USA}

\tableofcontents{}
\begin{abstract}
We develop an embedded boundary method (EBM) to solve the two-phase
incompressible flow with piecewise constant density. The front tracking
method is used to track the interface. The fractional step methods
are used to solve the incompressible Navier-Stokes equations while
the EBM is used in the projection step to solve an elliptic interface
problem for the pressure with a jump equal to the surface tension
force across the interface. Several examples are used to verify the
accuracy of the method.
\end{abstract}

\section{Introduction }

\label{sec:introduction}

We consider here the two-phase incompressible Navier-Stokes equations
(Zhijun Tan et al. \cite{TanLeLimKhoo2009}) \begin{eqnarray}
\rho(u_{t}+(u\cdot\nabla)u) & = & -\nabla p+\nabla\cdot\mu(\nabla u+(\nabla u)^{T})+F(x,t)+g(x,t),\label{eq:Navier-Stokes-Equations}\\
\nabla\cdot u & = & 0,\nonumber \end{eqnarray}
 where $\rho$ is the density (constant in each phase), $u$ is the
velocity field, $p$ is the pressure, $\mu$ is the viscosity, $g(x,t)$
is the external body force and \[
F(x,t)=\int_{\Gamma}f(s,t)\delta(x-X(s,t))ds\]
is the singular interface force, concentrated on the interface with
parametrization $X(s,t)$.

There are many different methods for solving the two-phase incompressible
flow that also track the phase interface. Some of the most popular
methods are the immersed boundary method \cite{Peskin2002}, the front
tracking method \cite{UnvTry92,TryggvasonBunner2001}, the level set
method \cite{SussmanSmerekaOsher1994,Sethian1999,OsherFedkiw2002},
and the volume-of-fluid method \cite{PuckettAlmgrenBellMarcusRider1997,PilliodPuckett2004}.
To deal with the singular source term $F(x,t)$ such as the surface
tension in equation (\ref{eq:Navier-Stokes-Equations}), these methods
simply use a discrete delta function to transform the source term
$F(x,t)$ defined only on the interface to the grid cells near the
interface, and then solve the Navier-Stokes equations as a one-phase
problem.

Recently, sharp interface methods to solve the two-phase flows directly
using the jump conditions due to the singular force term have received
attention. The ghost-fluid method \cite{GlimmMarchesinMcBryan1981,LiuFedkiwKang2000,KangFedkiwLiu2000}
has been used to solve the elliptic boundary value/elliptic interface
problem, and then used to solve the two-phase incompressible flow
to more accurately compute the solution satisfying the jump conditions
for the Navier-Stokes equations across the material interface. However,
since the ghost-fluid method uses approximate jump conditions, it
can only achieve at most first order accuracy for the elliptic interface
problem. Another popular method is the immersed interface method \cite{LeVequeLi1994,LiIto2006,LeeLeVeque2003,TanLeLimKhoo2009}
which has been applied to solve elliptic interface problems, parabolic
interface problems, and two-phase incompressible flow with discontinuous
viscosity. However, it seems that this method has not been used to
solve variable density two-phase incompressible flow. In addition
to the basic jump conditions, the immersed interface method requires
higher-order derivative jump conditions, which are often difficult
to be derived.

In this paper, we develop an embedded boundary method (EBM) \cite{JohansenColella98}
to solve the two-phase incompressible Navier-Stokes equations with
piecewise constant density. This method was originally proposed \cite{JohansenColella98}
to solve the elliptic boundary value problem on irregular domains.
It has since been used to solve the heat equation \cite{McCorquodaleColellaJohansen2001,Colella2006}
and the incompressible flow on a time-dependent domain \cite{MillerTrebotich2012}
with second order accuracy. Recently it has been extended to solve
the elliptic interface problem with second order accuracy in two and
three dimension \cite{WangSamulyakGuo2010,CrockettColellaGraves2010}
where the jump conditions for the potential and its flux are used
in the discretization of the elliptic equations using the finite volume
method.

The front tracking code \noun{FronTier} developed at Stony Brook University
has been successfully used for solving compressible flow \cite{GliGroLi98,DuFixGli05,BoLiuGlimmLi2011}.
It uses the ghost-fluid method \cite{GlimmMarchesinMcBryan1981,LiuFedkiwKang2000,KangFedkiwLiu2000}
to obtain sharp interface solution. In this paper, we solve the two-phase
incompressible flow equations using the embedded boundary method,
with the front tracking method used to track the material interface.
The jump conditions due to the surface tension force are accurately
solved in the projection step using the EBM while the jump conditions
due to the discontinuous viscosity are disregarded for simplicity.
The main contribution of this paper is the development of a sharp
interface EBM for the two phase incompressible flow. We also demonstrate
the extension of the method to the cylindrical coordinate. An extension
of the method has been used to solve large density ratio incompressible
multiphase magnetohydrodynamic flows \cite{GuoWangSamulyak2013}. 

The paper is organized as follows. In section \ref{sec:projection_method},
we review briefly the projection method used to solve the incompressible
flow without considering the interface. In section \ref{sec:EBM-for-the-Incompressible-Flow},
we describe the embedded boundary method used to solve the incompressible
flow with an internal interface. In section \ref{sec:examples}, we
show some examples to verify our methods. Finally, we give the conclusion.
In the appendices, we discuss the consequence of not considering the
jump condition due to discontinuous viscosity, the code structure,
and parallelization.

\section{The Projection Method}

\label{sec:projection_method}

In this section, we give the basic description of the projection methods
used to solve the incompressible flow without considering the interface.
Zhou et al. \cite{ZhoRayLim12} presented the verification of the
basic code for solving one phase flow and an implementation of the
two phase flow using the immersed boundary method with the front tracking
method.

To solve the Navier-Stokes equation, one of the most popular methods
is the projection method. There are many variations of the projection
method \cite{BrownCortezMinion2001}. In this paper, we use two different
but similar projection methods (denoted as PM1 and PM3 later on) with
the embedded boundary method to solve the two-phase incompressible
flow (PM1 is first order accurate and PM2 is second order accurate
when there is no interface). The projection methods for solving the
Navier-Stokes equations (\ref{eq:Navier-Stokes-Equations}) consists
of two steps \cite{BrownCortezMinion2001}. The first step solves
the convection and diffusion term to compute an intermediate velocity
at the new time step. The second step then solves the Poisson equation
to obtain the pressure, and use the new pressure to calculate the
new time step velocity, which satisfies the divergence free condition.
One salient feature of the two projection methods used in this paper
is that the pressure instead of the pressure increment is calculated
in the projection step. This feature makes it easier to use the EBM.
The reason is that for two phase incompressible flow with surface
tension, it is the pressure (not the pressure increment) which has
a jump equal to the surface tension across the interface.

\subsection{Projection Methods}

We describe three different variations of the projection methods.
However, the last projection method (denoted as PM3 later on) is presented
here only for comparison purpose and is not used in our algorithms.
For simplicity, we use $L$ to denote the discretization of the heat
operator in the Navier-Stokes equations.

\subsubsection{Projection Method PM1}

This projection method is similar to Chorin's original projection
method \cite{Chorin68,BrownCortezMinion2001}. \[
u^{*}=L(u^{n},-u\cdot\nabla u)\]

\begin{equation}
\frac{u^{n+1}-u^{*}}{\Delta t}=-\frac{1}{\rho}\nabla p^{n+\frac{1}{2}}\label{eq:PM1-ProjectionStep}\end{equation}
where $n$ is the time step index, $\Delta t$ is the time step size,
$u^{*}$ is the intermediate velocity, and $u^{n+1}$ is the velocity
at the new time step satisfying the divergence free property. Note
that this method is only first order accurate. A notable characteristics
of this method is that the projection step solves for the pressure
$p^{n+\frac{1}{2}}$ instead of the pressure increment.

\subsubsection{Projection Method PM2}

A second order accurate method \cite{MillerTrebotich2012} is the
following \begin{equation}
\tilde{u}=L(u^{n},-u\cdot\nabla u,-\nabla p^{n-\frac{1}{2}})\label{eq:PM2-HeatOperator}\end{equation}
\begin{equation}
\frac{u^{*}-\tilde{u}}{\Delta t}=\frac{1}{\rho}\nabla p^{n-\frac{1}{2}}\label{eq:PM2-PreProjectionStep}\end{equation}

\begin{equation}
\frac{u^{n+1}-u^{*}}{\Delta t}=-\frac{1}{\rho}\nabla p^{n+\frac{1}{2}}\label{eq:PM2-ProjectionStep}\end{equation}
This method has the same property as PM1 in that the projection step
solves for the pressure $p^{n+\frac{1}{2}}$ instead of the pressure
increment. This makes it easy to use our EBM method when a jump condition
about the pressure (instead of pressure increment) is given.

\subsubsection{Projection Method PM3}

Another popular second order accurate method is given by Bell, et
al \cite{BelColGlaz89,BrownCortezMinion2001}:\[
u^{*}=L(u^{n},-u\cdot\nabla u,-\nabla p^{n-\frac{1}{2}})\]

\begin{equation}
\frac{u^{n+1}-u^{*}}{\Delta t}=-\frac{1}{\rho}\nabla\phi^{n+\frac{1}{2}}\label{eq:eq:ProjectionMethod3-ProjectionStep}\end{equation}
\[
p^{n+\frac{1}{2}}=p^{n-\frac{1}{2}}+\phi^{n+\frac{1}{2}}\]
Note that the projection step solves for the pressure increment $\phi^{n+\frac{1}{2}}$.
For this reason, it is difficult to be used with a jump condition
for the pressure itself (such as the surface tension between two phases).
Thus, this method is not used in this paper.

\subsection{Time Discretization for the Calculation of the Intermediate Velocity}

There are many different methods to discretize the heat operator.
In the following, we briefly describe the Crank-Nicolson method and
the Additive Runge-Kutta Method that we have used in our method.

\subsubsection{Crank-Nicolson Method }

Using the Crank-Nicolson method for the heat operator and using PM1
as example, we have

\begin{equation}
\rho\left(\frac{u^{*}-u^{n}}{\Delta t}\right)=-(u\cdot\nabla)u^{n+\frac{1}{2}}+\frac{1}{2}(\mu\nabla^{2}u^{*}+\mu\nabla^{2}u^{n})+g\label{eq:projection_diffusion}\end{equation}
where the convection term is solved using an explicit Godunov-type
scheme as in \cite{BellColellaGlaz1989}: 
\begin{itemize}
\item extrapolate the velocity to the cell face at time $n+\frac{1}{2}$
(the following shows the formula for only one face): \begin{equation}
u_{i+1/2,j}^{n+1/2}=u_{i,j}^{n}+\frac{1}{2}\Delta x\left.\frac{\partial u}{\partial x}\right|_{i,j}^{n}+\frac{1}{2}\Delta t\left.\frac{\partial u}{\partial t}\right|_{i,j}^{n},\label{eq:convection-velocity-extrapolation}\end{equation}
 where $i$, $j$ denote the cell index, $(i+1/2,j$) denotes one
of the the cell face, $n$ denote the time step index, and $\frac{\partial u}{\partial t}$
can be replaced by the Navier-Stokes equations\[
u_{t}=-(u\cdot\nabla)u-\frac{1}{\rho}\nabla p+\mu(\nabla^{2}u)+g(x,t)).\]

\item solve the Riemann problem (Burgers' equation) to find the cell face
velocity at time $n+\frac{1}{2}$. 
\item use the cell face velocity to calculate the convection term $(u\cdot\nabla u)^{n+\frac{1}{2}}$. 
\end{itemize}

\subsubsection{Additive Runge-Kutta Method}

The implicit Runge-Kutta method (\cite{TwizellGumelArigu1996}) was
used to solve the time dependent parabolic initial boundary value
problem in (\cite{McCorquodaleColellaJohansen2001,Colella2006}) and
the parabolic interface problem in \cite{WangSamulyakGuo2010}. In
this paper, we instead use the additive Runge-Kutta method \cite{LiuZou2006}
which seems to be easier for the discretization of the Navier-Stokes
equation.

Using the notation in \cite{LiuZou2006}, we want to solve the following
ordinary differential equation (ODE):\begin{equation}
y'(t)=f(t,y)+g(t,y)\label{eq:additive-ODE}\end{equation}
 where $f$ is linear operator of $y$ and $g$ is a nonlinear operator.
We use an implicit scheme for $f$ and an explicit scheme for $g$.
When used with PM2, $f=\Delta u$ and $g=-u\cdot\nabla u-\nabla p^{n-\frac{1}{2}}$.
Note that to solve the convection term $u\cdot\nabla u$ for the additive
Runge-Kutta method, there is no need to do time extrapolation in equation
(\ref{eq:convection-velocity-extrapolation}). 

To solve the ODE (\ref{eq:additive-ODE}), the $s$-stage Runge-Kutta
method has the following form\[
y_{i}^{(n)}=y_{s}^{(n-1)}+h\sum_{j=1}^{s}a_{ij}f(t_{n-1}+c_{j}h,y_{j}^{(n)})+h\sum_{j=1}^{s}b_{ij}g(t_{n-1}+c_{j}h,y_{j}^{(n)})\]
 where $i=1,2,...,s$, $n=1,2,3,...$, and $c_{i}=\sum_{j=1}^{s}a_{ij}=\sum_{j=1}^{s}b_{ij}$.
The coefficients are generally written as the tableau in Table \ref{tab:Additive-Runge-Kutta-Coefficients}.

\begin{table}
\caption{\label{tab:Additive-Runge-Kutta-Coefficients}Additive Runge-Kutta
coefficients tableau}

\centering{}\begin{tabular}{|c|c|c|c|c|c|c|}
\hline 
$c_{1}$  & $a_{11}$  & ...  & $a_{1s}$  & $b_{11}$  & ...  & $b_{1s}$\tabularnewline
\hline
\hline 
...  & ...  &  & ...  & ...  &  & ...\tabularnewline
\hline 
$c_{s}$  & $a_{s1}$  & ...  & $a_{ss}$  & $b_{s1}$  & \multicolumn{1}{c||}{... } & $b_{ss}$\tabularnewline
\hline
\end{tabular}
\end{table}

In this paper, we use the $L$-stable two stage additive Runge-Kutta
scheme shown in Table \ref{tab:L-stable-two-stage-Runge-Kutta}.

\begin{table}
\caption{\label{tab:L-stable-two-stage-Runge-Kutta}L-stable two stage additive
scheme}

\centering{}\begin{tabular}{|c|c|c|c|c|c|c|}
\hline 
0  & $0$  & 0  & 0  & 0  & 0  & 0\tabularnewline
\hline
\hline 
$\frac{1}{2}$  & $\frac{1}{2}(-1+\sqrt{2})$  & $1-\frac{\sqrt{2}}{2}$  & 0  & $\frac{1}{2}$  & 0  & 0\tabularnewline
\hline 
1  & $1-\frac{\sqrt{2}}{2}$  & $\sqrt{2}$-1  & $1-\frac{\sqrt{2}}{2}$  & 0  & 1  & 0\tabularnewline
\hline
\end{tabular}
\end{table}

\subsection{Calculation of the New Velocity}

Now we consider the projection step of PM1 and PM2 for calculating
the pressure and the new divergence free velocity. For a regular grid
cell containing no interface, the new velocity is calculated using:

\begin{equation}
\rho\left(\frac{u^{n+1}-u^{*}}{\Delta t}\right)=-\nabla p^{n+1/2}\label{eq:projection}\end{equation}
 where $p^{n+1/2}$ denotes the new pressure and is calculated using

\begin{equation}
\nabla\cdot\left(\frac{1}{\rho}\nabla p^{n+1/2}\right)=\frac{1}{\Delta t}\nabla\cdot u^{*},\label{eq:projection_elliptic}\end{equation}
 and the boundary condition

\begin{equation}
n\cdot\nabla p^{n+1/2}=0.\label{eq:boundaryCondition_dp}\end{equation}

We have used $\nabla\cdot u^{n+1}=0$ and equation (\ref{eq:projection})
to obtain the Poisson equation (\ref{eq:projection_elliptic}).

\section{The EBM for Two-Phase Incompressible Flow}

\label{sec:EBM-for-the-Incompressible-Flow}

In this section, we first review the jump conditions of the incompressible
flow with an internal interface. Next we review our embedded boundary
method for solving the elliptic interface problem in subsection \ref{sub:review-of-embedded-boundary-method}.
Then we show how to use the EBM in a cylindrical coordinate system.
The method to calculate the pressure gradient used in PM2 is described
in subsection \ref{sub:Calculation-of-Pressure-Gradient}. At last
we apply the EBM to solve the elliptic interface problem of the projection
step for the two-phase incompressible flow in subsection \ref{sub:EBM-for-projection-step}.

\subsection{The Jump Conditions across an Internal Interface}

\label{sub:the_jump_conditions}

For two phase flow with discontinuous coefficients, the solutions
satisfy the following jump conditions across the internal interface:

\begin{equation}
[u]=0\label{eq:jump_u}\end{equation}

\begin{equation}
[p]=2[\mu\frac{\partial u}{\partial n}]\cdot n+f\cdot n\label{eq:jump_p}\end{equation}

\begin{equation}
[\mu\frac{\partial u}{\partial n}]\cdot\tau+[\mu\frac{\partial u}{\partial\tau}]\cdot n+f\cdot\tau=0\label{eq:jump_du1}\end{equation}

\begin{equation}
[\frac{\partial u}{\partial n}]\cdot n=0\label{eq:jump_du2}\end{equation}

\begin{equation}
[\frac{1}{\rho}\frac{\partial p}{\partial n}]=0\label{eq:jump_dp}\end{equation}
 where $n$ is the normal vector and $\tau$ is the tangential vector.
Note that we have four jump conditions for the two velocity components
in $2D$ (six in $3D$) and two jump conditions for the pressure. 
\begin{itemize}
\item Jump condition (\ref{eq:jump_u}) is due to the viscous flow (Zhi,
et al \cite{TanLeLimKhoo2009}). 
\item Jump conditions (\ref{eq:jump_p}, \ref{eq:jump_du1}) are the result
of force balancing in the normal and tangential directions (Ito and
Li \cite{ItoLi2006}). The derivation of the jump conditions (\ref{eq:jump_p},\ref{eq:jump_du1})
is given by Ito and Li \cite{ItoLi2006} for the Stokes equations.
It is trivial to extend it for Navier-Stokes equation. 
\item Jump condition (\ref{eq:jump_du2}) is due to $[\nabla\cdot u]=0$
(Lai and Li \cite{LaiLi2001}). 
\item Jump condition (\ref{eq:jump_dp}) is needed to obtain a continuous
velocity in the normal direction across the interface when solving
the elliptic interface problem for the pressure.
\end{itemize}
In this paper, we assume that the singular force term on the interface
has only a normal component, which means that $f\cdot\tau=0$. The
surface tension force satisfies this assumption. For simplicity, we
also ignore the jump conditions due to the discontinuous viscosity.
Appendix \ref{sub:Appendix-Approximate-Jump-Conditions} shows the
consequence of such a simplification. Therefore the velocity and its
derivatives are continuous across the interface, while the pressure
has a jump equal to the surface tension across the interface.

\subsection{The Embedded Boundary Method for the Elliptic Interface Problem}

\label{sub:review-of-embedded-boundary-method}

We review briefly the embedded boundary method for solving the elliptic
interface problem. For more detail, refer to \cite{WangSamulyakGuo2010,CrockettColellaGraves2010,JohansenColella98,McCorquodaleColellaJohansen2001}.
We assume that the interface can only cross any cell edge at most
once, as implicitly assumed in the Marching Cubes algorithm \cite{LorCli87}.
This assumption is necessary to limit the number of different cases
that could arise. The geometric information needed of the partial
cells needed by the embedded boundary method is calculated using the
divergence theorem for each cell case by case \cite{Wang2012}.

The elliptic interface problem is a special elliptic problem with
an internal interface: \begin{equation}
\nabla\centerdot\frac{\nabla p}{\rho}=f\label{eq:elliptic}\end{equation}
 where $\rho$ is a piecewise continuous function with jump across
the internal interface and $f$ is a given function which is continuous
inside each part of the domain. To close the problem, boundary conditions
are needed for the exterior and interior boundary. Either Dirichlet
or Neumann boundary can be given on the exterior boundary. For the
interior interface, we have the following two jump conditions: \begin{equation}
[p]=J_{1}(\mathbf{x}),\label{eq:elliptic_potential_jump}\end{equation}
 \begin{equation}
[\frac{1}{\rho}\frac{\partial p}{\partial n}]=J_{2}(\mathbf{x}),\label{eq:elliptic_flux_jump}\end{equation}
 where $J_{1}$ and $J_{2}$ are given functions of the spatial variables
\cite{LeVequeLi1994}. Note that our method could be easily extended
to more general cases where $J_{1}$ and $J_{2}$ are functions of
the unknown variables $p$ defined on the interface.

The EBM can be used to solve the elliptic problem with an irregular
domain boundary and an internal interface. It uses a Cartesian mesh.
The mesh cells are classified into four types: 
\begin{itemize}
\item An \emph{external} cell is outside of the computational domain and
thus is not used in the computation. 
\item An \emph{internal} cell is a cell wholly located inside the computational
domain, possibly with one of its cell face being the domain boundary. 
\item A \emph{boundary} cell is a cell intersected by the irregular exterior
domain boundary with part of the cell out of the domain and part of
the cell inside the domain. 
\item A \emph{partial} cell is a cell intersected by the internal interface.
It is separated into two or more parts by the internal interface.
Those different parts are also called partial cells in the following.
\end{itemize}
\begin{table}
\begin{centering}
\caption{\label{tab:Number-of-Cell-Unknowns}Number of material components
and cell unknowns for different cell types}

\par\end{centering}

\centering{}\begin{tabular}{|c|c|c|c|}
\hline 
 & \multicolumn{3}{c|}{Number of}\tabularnewline
\hline 
\multicolumn{1}{|c|}{Cell Type} & Components & Center Unknowns & Interface Unknowns\tabularnewline
\hline
\hline 
External  & N/A  & N/A  & N/A\tabularnewline
\hline 
Internal  & 1  & 1  & 0\tabularnewline
\hline 
Boundary  & 1  & 1  & 0\tabularnewline
\hline 
Partial  & 2  & 2  & 2\tabularnewline
\hline
\end{tabular}
\end{table}

\begin{figure}
\centering{ \includegraphics[scale=0.4]{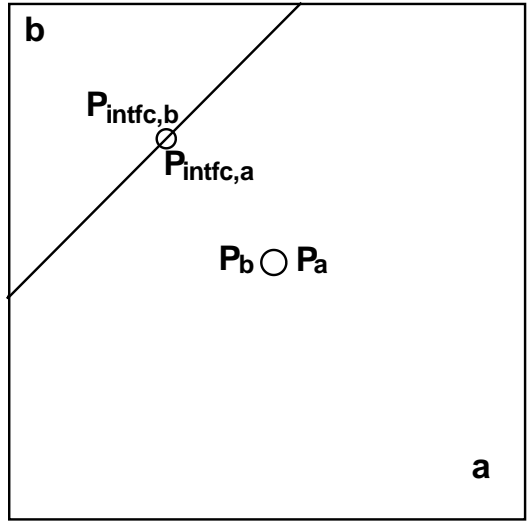}

\caption{\label{fig:the-cell-unknowns-for-ebm}Placement of unknowns in a partial
cell containing the cell internal interface. $p_{a}$ and $p_{b}$
are cell center unknowns for component $a$ and $b$ respectively.
Similarly, $p_{intfc,a}$ and $p_{intf,b}$ are cell interface unknowns
for component $a$ and $b$ respectively.}

} 
\end{figure}

When the EBM is used to solve the elliptic interface problem, one
or more unknowns are defined at the cell center, as shown in Table
\ref{tab:Number-of-Cell-Unknowns}. For an interior cell or a boundary
cell, only one unknown is needed at the cell center and a standard
finite volume method can be used to setup one algebraic equation using
the elliptic equation as in \cite{WangSamulyakGuo2010,JohansenColella98,McCorquodaleColellaJohansen2001,Colella2006}.
For a partial cell, four unknowns in total are needed to make the
discretization of the elliptic equation consistent with the two interface
jump conditions (\ref{eq:elliptic_potential_jump}), (\ref{eq:elliptic_flux_jump}).
Figure \ref{fig:the-cell-unknowns-for-ebm} depicts the placement
of unknowns in a partial cell with internal interface. The cell contains
two partial cells (for material components $a$ and $b$) which are
separated by the interface. For each partial cell, there is one unknown
defined at the cell center ($p_{a}$ for component $a$, $p_{b}$
for component $b$) for the discretization of the elliptic equation.
In order to satisfy the two jump conditions (\ref{eq:elliptic_potential_jump},
\ref{eq:elliptic_flux_jump}), two more unknowns $p_{intfc,a}$ and
$p_{intfc,b}$ are defined at the center of the cell internal interface
(portion of the interface contained within the cell) for the two components
$a$ and $b$. Thus, four unknowns are defined in total for one partial
cell. Four algebraic equations can be constructed using the elliptic
equation for the two partial cells and the two jump conditions across
the cell internal interface 

\subsubsection{Discretization of the Jump Conditions}

\label{sub:Discretization-Of-The-Jump-Conditions}

\begin{figure}
\centering{ \includegraphics[scale=0.3]{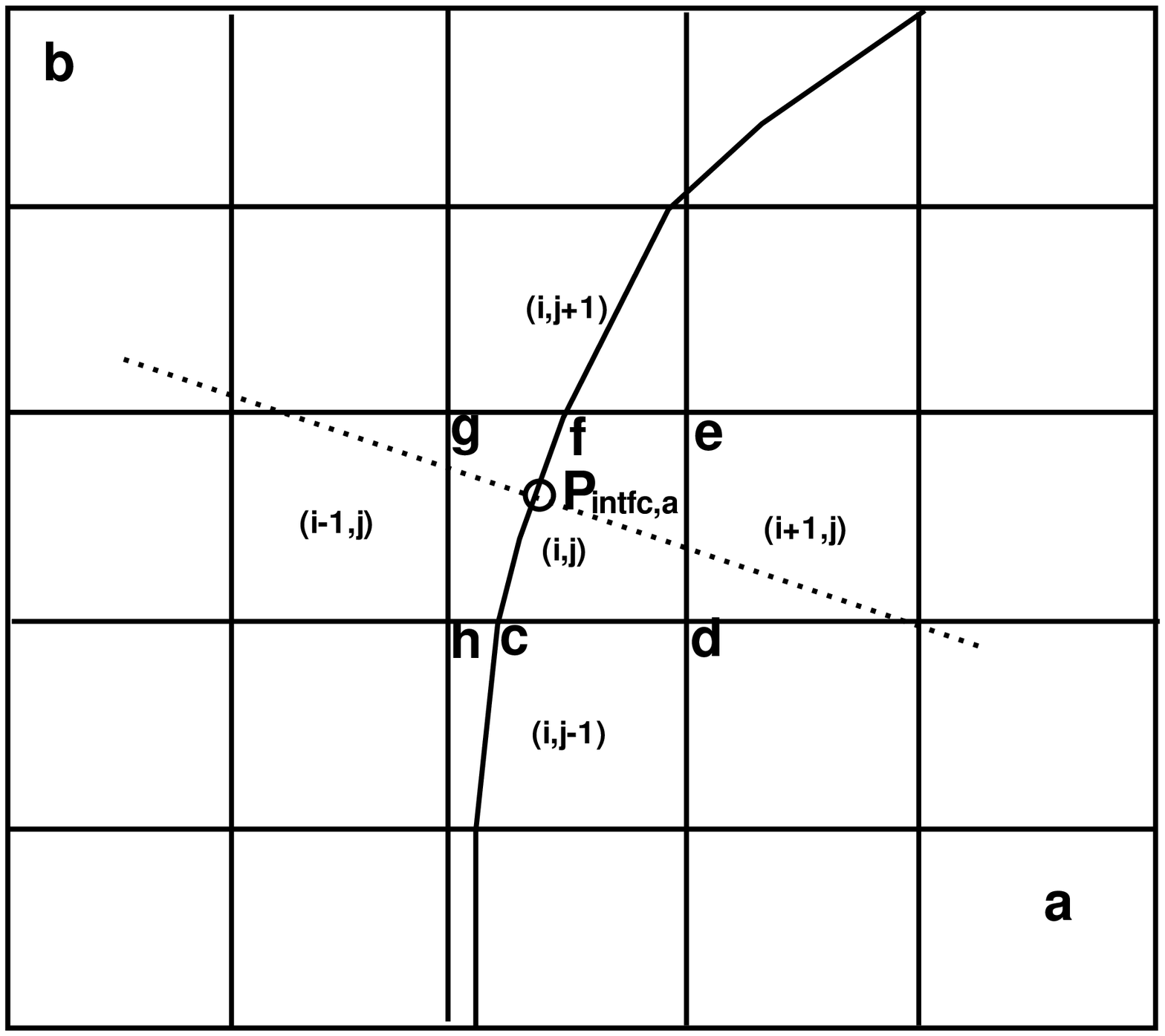}

\caption{\label{fig:stencil-for-the-intfc-unknowns}Stencil for the discretization
for a partial cell}

} 
\end{figure}

We first describe the method for the discretization of the two jump
conditions across the interface for the cell $(i,j)$. A schematic
of the corresponding stencil and states used in the interpolation
method is shown in Figure \ref{fig:stencil-for-the-intfc-unknowns}.
Two unknowns $p_{intfc,a}$ and $p_{intfc,b}$ are defined at the
center of the cell interface for components $a$ and $b$ respectively.

We assume the direction of the normal to the interface as pointing
from $a$ to $b$. For the first jump condition (\ref{eq:elliptic_potential_jump}),
the discretization is simply \begin{equation}
p_{intfc,b}-p_{intfc,a}=J_{1}\label{eq:elliptic_discretized_equation_1}\end{equation}
 To discretize the second jump condition (\ref{eq:elliptic_flux_jump}),
we need to calculate the normal direction derivatives of the unknowns.
There are many different approach available as in \cite{WangSamulyakGuo2010,JohansenColella98,McCorquodaleColellaJohansen2001}.
The main idea is to construct a polynomial using the unknowns with
the same component and then take derivative along the normal direction
to get the flux for that component. To construct the polynomial for
component $a$, we pick cell interface unknown $p_{intfc,a}$ from
cell $(i,j)$ and then cell center unknowns from other neighboring
cells. In $2D$, we need $3$ unknowns to construct a linear polynomial
(for first order accurate flux) and $6$ unknowns to construct a quadratic
polynomial (for second order accuracy flux). Least square fitting
is also possible.

Taking the normal derivative of the constructed polynomials, we obtain
the normal derivatives at the cell interface center for the component
$a$, $\left.\frac{\partial p}{\partial n}\right|_{a}$, and the component
$b$, $\left.\frac{\partial p}{\partial n}\right|_{b}$. Using the
second jump condition (\ref{eq:elliptic_flux_jump}), we have \begin{equation}
\left.\frac{1}{\rho}\right|_{b}\cdot\left.\frac{\partial p}{\partial n}\right|_{b}-\left.\frac{1}{\rho}\right|_{a}\cdot\left.\frac{\partial p}{\partial n}\right|_{a}=J_{2}\label{eq:elliptic_discretized_equation_2}\end{equation}

\paragraph{Method for calculating the flux across the cell internal interface}

We use linear polynomial construction in $2D$ as an example. A linear
polynomial in $2D$ can be written as \[
p(x,y)=\beta_{0}+\beta_{1}x+\beta_{2}y\]
 where $\beta_{0}$, $\beta_{1}$ and $\beta_{2}$ are undetermined
coefficients. Thus we need to find $3$ cell unknowns ($1$ unknown
should be a cell interface unknown) for constructing the polynomial.
Denoting the coordinate locations of the unknowns as $(x_{1},y_{1})$,
$(x_{2},y_{2})$, and $(x_{3},y_{3})$, the unknowns as $p_{1}$,
$p_{2}$, and $p_{3}$, we have \[
\left(\begin{array}{ccc}
1 & x_{1} & y_{1}\\
1 & x_{2} & y_{2}\\
1 & x_{3} & y_{3}\end{array}\right)\left(\begin{array}{c}
\beta_{0}\\
\beta_{1}\\
\beta_{2}\end{array}\right)=\left(\begin{array}{c}
p_{1}\\
p_{2}\\
p_{3}\end{array}\right).\]

Denoting the above matrix as $X$, the coefficient vector as $\beta$,
and the right hand side as $P$, we have \[
X\beta=P,\]
 or\[
\beta=X^{-1}P.\]

Now the constructed polynomial can be written as \[
\left(\begin{array}{ccc}
1 & x & y\end{array}\right)\left(\begin{array}{c}
\beta_{0}\\
\beta_{1}\\
\beta_{2}\end{array}\right).\]
 The derivatives of the polynomial can be written as\begin{eqnarray}
p_{x} & = & \beta_{1}\nonumber \\
p_{y} & = & \beta_{2}\nonumber \\
\frac{\partial p}{\partial n} & = & p_{x}n_{1}+p_{y}n_{2}=\left(\begin{array}{ccc}
0 & n_{1} & n_{2}\end{array}\right)\left(\begin{array}{c}
\beta_{0}\\
\beta_{1}\\
\beta_{2}\end{array}\right)=\hat{n}^{T}\beta\label{eq:directional-directive}\end{eqnarray}
 where $n=(\begin{array}{cc}
n_{1}, & n_{2}\end{array})^{T}$ is the unit normal vector, and $\hat{n}=\left(\begin{array}{ccc}
0, & n_{1}, & n_{2}\end{array}\right)^{T}$. Thus we have \begin{eqnarray*}
\frac{\partial p}{\partial n} & = & \hat{n}^{T}X^{-1}P\\
 & = & ((X^{-1})^{T}\hat{n})^{T}P,\\
 & = & ((X^{T})^{-1}\hat{n})^{T}P\end{eqnarray*}
 Therefore, we can represent $\frac{\partial p}{\partial n}$ as a
linear combination of the unknowns $P$. Similarly, we can calculate
the directional derivative using second degree polynomial construction.

Now suppose that instead we want to use least square fitting, we have
\[
\beta=((X^{T}X)^{-1}X^{T})P.\]
 The directional derivative can be calculated as \begin{eqnarray*}
\frac{\partial p}{\partial n} & = & \hat{n}^{T}\beta=\hat{n}^{T}((X^{T}X)^{-1}X^{T})P\\
 & = & (X(X^{T}X)^{-1}\hat{n})^{T}P,\end{eqnarray*}
 thus $\frac{\partial p}{\partial n}$ can also be written as a linear
combination of the unknowns $P$.

\subsubsection{Equations for the Cell Center Unknowns}

For the two cell center unknowns defined at the partial cell $(i,j)$,
we can use the EBM technique to set up two algebraic equations by
integrating the elliptic equation over the corresponding partial cells.
See Figure \ref{fig:stencil-for-the-intfc-unknowns} for the stencil
to set up the equation for the unknown of the component $a$ using
the partial cell $cdef$.

Integrating the equation (\ref{eq:elliptic}) over the partial cell
$cdef$ and using the divergence theorem, we obtain the following
expressions: \[
\iint_{cdef}\nabla\cdot\frac{\nabla p}{\rho}dxdy=\oint_{\partial(cdef)}\frac{\nabla p}{\rho}\cdot nds=\iint_{cdef}fdxdy,\]
 or \[
\int_{cd}\frac{\nabla p}{\rho}\cdot nds+\int_{de}\frac{\nabla p}{\rho}\cdot nds+\int_{ef}\frac{\nabla p}{\rho}\cdot nds+\int_{fc}\frac{\nabla p}{\rho}\cdot nds=\iint_{cdef}fdxdy,\]
 which is \begin{equation}
l_{cd}\cdot flux_{cd}+l_{de}\cdot flux_{de}+l_{ef}\cdot flux_{ef}+l_{fc}\cdot flux_{fc}=\iint_{cdef}fdxdy\end{equation}
 where $l_{mn}$ is the length between $m$ and $n$. Therefore we
only need to calculate the flux across the cell edges or cell interface.
For $flux_{cd}$, a second order derivative is calculated by using
a linear interpolation of $\frac{p_{i,j-1}-p_{i,j}}{\triangle x}$
and $\frac{p_{i+1,j-1}-p_{i+1,j}}{\triangle x}$ to the center of
$cd$ (see \cite{WangSamulyakGuo2010,JohansenColella98}). For $flux_{de}$,
we simply use $\frac{p_{i+1,j}-p_{i,j}}{\triangle y}$ to calculate
the derivative. For $flux_{ef}$, a linear interpolation of $\frac{p_{i,j+1}-p_{i,j}}{\triangle x}$
and $\frac{p_{i+1,j+1}-p_{i+1,j}}{\triangle x}$ to the center of
$ef$ is used. And $flux_{fc}$ is calculated by $\frac{\partial p}{\partial n}|_{a}$
used in (\ref{eq:elliptic_discretized_equation_2}). Note that we
need to multiply the derivatives calculated with $\frac{1}{\rho}$
at the corresponding edge or cell interface center to obtain the flux.
In the same way, we calculate fluxes for the other partial cells.
More details can be found in \cite{WangSamulyakGuo2010,JohansenColella98}.

\subsubsection{Geometric Property Calculation}

We use $3D$ as example. Since the EBM is a finite volume method,
we need to calculate the partial cell volumes accurately. To solve
the elliptic interface problem, we also need to calculate the cell
interface area, normal and center.

Due to the assumption that the interface crosses the cell edge at
most once, there are finite number of cases for the partial cell configurations
with two components. The cell volume, cell interface area, interface
normal, and interface center are calculated case by case by using
the divergence theorem:

\begin{equation}
\int_{\Omega}\nabla\cdot Fdv=\int_{\partial\Omega}\overrightarrow{n}\cdot Fds.\label{eq:Divergence-Theorem}\end{equation}
 where $\overrightarrow{F}=(F_{x},F_{y},F_{z})$. Thus, we can use
a cell boundary integration instead of a cell volume integration to
calculate the geometric information needed.

In the Cartesian coordinate, the divergence operator is defined as
\[
\nabla\cdot\overrightarrow{F}=\frac{\partial F_{x}}{\partial x}+\frac{\partial F_{y}}{\partial y}+\frac{\partial F_{z}}{\partial z}.\]

To calculate the volume of the domain $\Omega$, we can let $\overrightarrow{F}=\overrightarrow{X}=(x,y,z)^{T}$
where $\overrightarrow{X}$ is the Cartesian coordinate. Note that
$F$ is not unique. Then we have \begin{eqnarray*}
\int_{\Omega}\nabla\cdot\overrightarrow{F}dv & = & \int_{\partial\Omega}n\cdot\overrightarrow{F}ds\\
\int_{\Omega}3dv & = & \int_{\partial\Omega}n\cdot\overrightarrow{X}ds\end{eqnarray*}
 Therefore, we can calculate the volume using the surface integration\begin{equation}
Volume=\int_{\Omega}dv=\frac{\int_{\partial\Omega}n\cdot\overrightarrow{X}ds}{3}.\label{eq:volume-calculation-for-Cartesian-Coordinate}\end{equation}
 For the cell interface area, normal and center, we do surface integrations:
\[
\int_{intfc}ds\]
 \[
\int_{intfc}\overrightarrow{n}ds\]
 \[
\int_{intfc}\overrightarrow{X}ds\]
 where $intfc$ refers to the cell interface. For more detail, refer
to \cite{Wang2012}.

\subsection{The Embedded Boundary Method in Cylindrical Coordinates}

\label{sec:EBM_in_cylindrical_coordinate}

In this subsection, we briefly describe the embedded boundary method
used to solve an elliptic interface problem in a cylindrical coordinate.

To discretize the elliptic equation (\ref{eq:elliptic}) for a mesh
cell, we use the divergence theorem to get \begin{eqnarray*}
\oint\frac{\nabla p}{\rho}\cdot nds & = & \int fdv\end{eqnarray*}
 or\[
\oint\frac{1}{\rho}\frac{\partial p}{\partial n}ds=\int fdv\]

This is true for all coordinate systems. However, we need to change
the formula for the calculation of the geometric property for the
mesh cell, such as the cell face area, cell volume, cell interface
area, normal and the center. The formula for the flux calculations
across the cell face and the cell internal interface also need to
be modified.

\subsubsection{Finite Difference Scheme for the Flux Calculation}

To calculate the flux across the cell faces or partial cell interface,
we use \[
\frac{\partial p}{\partial n}=\nabla p\cdot\overrightarrow{n}\]
 where the gradient operator is defined as \[
\nabla=\overrightarrow{r}\frac{\partial}{\partial r}+\overrightarrow{\theta}\frac{1}{r}\frac{\partial}{\partial\theta}+\overrightarrow{z}\frac{\partial}{\partial z}\]
 and $\overrightarrow{n}$ is the normal of the cell faces.

For the flux across the partial cell interface, we need to modify
the equation \ref{eq:directional-directive} using the Cylindrical
coordinate directional derivative instead.

\subsubsection{Geometric Property Calculation}

The surface elements in cylindrical coordinate for the cell face of
the internal cell are \[
ds=rd\theta dz,\]
\[
ds=drdz,\]
\[
ds=rdrd\theta\]
 in the respective coordinate planes. They can be used to calculate
the cell face area for the $6$ cell faces. The volume element is
\[
dv=rdrd\theta dz.\]

To calculate the partial cell volume, we also use the divergence theorem.
In a cylindrical coordinate, the divergence operator is defined as
\[
\nabla\cdot\overrightarrow{F}=\frac{1}{r}\frac{\partial}{\partial r}(rF_{r})+\frac{1}{r}\frac{\partial F_{\theta}}{\partial\theta}+\frac{\partial F_{z}}{\partial z}.\]

To calculate the volume of the domain $\Omega$, we can let $\overrightarrow{F}=(r,r\theta,z)^{T}$,
and then we have \begin{eqnarray*}
\int_{\Omega}\nabla\cdot\overrightarrow{F}dv & = & \int_{\partial\Omega}n\cdot\overrightarrow{F}ds\\
\int_{\Omega}4dv & = & \int_{\partial\Omega}n\cdot\overrightarrow{F}ds\end{eqnarray*}
Note that the choice for $\overrightarrow{F}$ is not unique. Therefore,
we can calculate the volume using the surface integration\begin{equation}
Volume=\int_{\Omega}dv=\frac{\int_{\partial\Omega}n\cdot\overrightarrow{F}ds}{4}.\label{eq:volume-calculation-for-Cylindrical-Coordinate}\end{equation}

For more detail, refer to \cite{Wang2012}.

\subsection{Calculation of the Pressure Gradient}

\label{sub:Calculation-of-Pressure-Gradient}

To calculate the pressure gradient used in (\ref{eq:PM2-HeatOperator})
and (\ref{eq:PM2-PreProjectionStep}), we simply use least square
fitting. To calculate the gradient for cell $(i,j)$ with component
$c$, we do the following:
\begin{enumerate}
\item find the nearest cell $(\hat{i},\hat{j})$ containing component $c$.
\item collect all the states with the same component from neighboring cells
with cell index $(\tilde{i},\tilde{j})$ satisfying $\left|\tilde{i}-\hat{i}\right|\leq2$
and $\left|\tilde{j}-\hat{j}\right|\leq2$.
\item use least square fitting to fit a quadratic polynomial and take gradient
to calculate the pressure gradient needed. See subsection \ref{sub:Discretization-Of-The-Jump-Conditions}
for detail.
\end{enumerate}

\subsection{The EBM for the Projection Step}

\label{sub:EBM-for-projection-step}

We solve the elliptic equation (\ref{eq:projection_elliptic}) in
the projection step. The EBM for a general elliptic interface problem
is given in previous subsections. The two jump conditions given by
(\ref{eq:jump_p}, \ref{eq:jump_dp}) are used to connect the pressure
solution across the interface.

Thus, for the projection step, the equation is \begin{equation}
\nabla\cdot(\frac{1}{\rho}\nabla(p^{n+1/2}))=\frac{1}{\Delta t}\nabla\cdot u^{*}\end{equation}
 coupled with the two jump conditions (omitting the effect due to
the viscosity jump across the interface) \begin{equation}
[p^{n+1/2}]=f\cdot n\end{equation}
and\begin{equation}
[\frac{1}{\rho}\frac{\partial p^{n+1/2}}{\partial n}]=0\end{equation}
where $f\cdot n=\sigma\kappa$ for the surface tension force between
the two phases ($\sigma$ is the surface tension coefficient and $\kappa$
is the mean curvature).

Our projection method consists of the following steps: 
\begin{enumerate}
\item Calculate cell face velocity $u_{face}^{*}$ using intermediate cell
center velocity $u^{*}$ by linear interpolation. 
\item Calculate the elliptic interface equation for $p$:\[
\begin{cases}
\nabla\cdot\frac{\nabla p}{\rho}=\frac{1}{\Delta t}\nabla\cdot u_{face}^{*}\\
{}[p]=\sigma\kappa\\
{}[\frac{1}{\rho}\frac{\partial p}{\partial n}]=0\end{cases}\]
where the right hand side of the elliptic equation is calculated using
the divergence theorem.
\item Calculate the projected cell face velocity $u_{face,c}^{n+1}$ for
each full/partial cell face using:\[
u_{face,c}^{n+1}=u_{face}^{*}-\frac{\Delta t}{\rho_{face,c}}\cdot\left.\frac{\partial p}{\partial n}\right|_{face,c}\]
where $\nabla p_{face,c}$ is calculated using the cell center pressure
for component $c$.
\item Calculate the averaged cell face velocity $u_{face}^{n+1}$ using
$u_{face,c}^{n+1}$. For example, in Figure \ref{fig:stencil-for-the-intfc-unknowns},
\[
u_{ge}^{n+1}=\frac{l_{gf}u_{gf,b}^{n+1}+l_{fe}u_{fe,a}^{n+1}}{l_{ge}}\]
 where $l_{ge}$, $l_{fg}$, $l_{fe}$ are the length between corresponding
point on the cell face, $u_{gf,b}^{n+1}$ is the velocity through
$gf$ and $u_{fe,a}^{n+1}$ is the velocity through $fe$. 
\item Calculate the cell centered velocity using $u_{face}^{n+1}$ by using
the $2nd$ order TVD reconstruction algorithms \cite{Balsara2001}.
A simple alternative is to use interpolation of the cell face velocity
to calculate the cell center velocity. 
\end{enumerate}

\section{Examples}

\label{sec:examples}

We compare theoretical and simulated bubble oscillation frequencies
in $2D$ and $3D$ to verify our embedded boundary method to solve
two-phase incompressible flow. To apply EBM to cylindrical coordinates,
we first verify the accuracy of the EBM by solving an elliptic interface
problem with known solution. Then we apply our methods to solve a
more complicated problem of engineering interest.

\subsection{Bubble Oscillation in Two and Three Dimension}

We consider the oscillation period of a droplet under zero gravity,
for which there is an analytical solution and the surface tension
is the dominant force. We have used both PM1 and PM2 to obtain the
simulation result. However, it is found that the difference between
the results obtained using the two methods are negligible.

\subsubsection{Bubble Oscillation in Two Dimension}

For a $2D$ droplet under zero gravity, when the initial position
of the droplet interface with small perturbation is given by \[
R(\theta)=R_{0}+\epsilon cos(n\theta)\]
where $R_{0}$ the unperturbed radius, $\epsilon$ is the amplitude
of the perturbation and $n$ is the order of the Legendre polynomial,
the oscillation frequency is given in \cite{FyfeOranFritts1988,ShiJur02}
as \[
\omega_{n}=\sqrt{\frac{(n^{3}-n)\sigma}{(\rho_{d}+\rho_{o})R_{0}^{3}}}\]
where $\rho_{d}$ and $\rho_{o}$ are the densities for the droplet
and outer fluid and $\sigma$ is the surface tension coefficient. 

\begin{figure}
\caption{\label{fig:Bubble-Oscillation-2D-Interface-Velocity}Bubble Oscillation
in $2D$, with interface velocity drawn as vector starting from the
interface at time $t=0.5$}

\centering{}\includegraphics[scale=0.3]{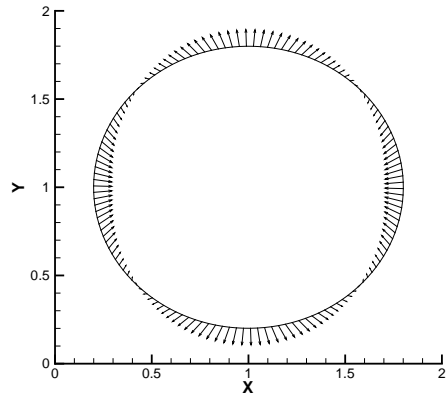}
\end{figure}

\begin{figure}
\caption{\label{fig:Bubble-Oscillation-2D-Pressure}Bubble Oscillation in $2D$,
Pressure at time $t=0.5$}

\centering{}\includegraphics[scale=0.3]{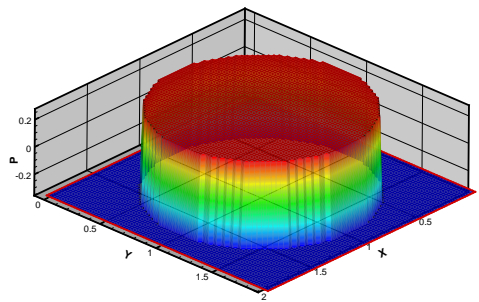} 
\end{figure}

The period is given by \[
T_{n}=\frac{2\pi}{\omega_{n}}.\]
We run the simulation with the domain as $[0,2]\times[0,2]$, the
densities and dynamic viscosities of the droplet and the outer fluid
given by $\rho=\{1,0.05\}$, $\nu=\{0.0005,0.0005\times0.01\}$. We
impose a surface tension coefficient of $\sigma=0.5$. For the initial
interfacial position, we use $R_{0}=0.8$, $n=2$ and $\epsilon=0.05$.
Table \ref{tab:Bubble-Oscillation-2D-1} shows the convergence of
the oscillation period under mesh refinement. Figure \ref{fig:Bubble-Oscillation-2D-Interface-Velocity}
shows the interface velocity of the droplet at time $t=0.5$. From
the picture we can see that the interface velocity changes smoothly
along the interface. Figure \ref{fig:Bubble-Oscillation-2D-Pressure}
shows the pressure over the whole computational domain. It is apparent
that the pressure has a jump across the interface due to the surface
tension. Due to the small perturbation of the initial interface ($\epsilon=0.05$)
and large radius of the droplet ($R_{0}=0.8$), the variation of the
surface tension along the interface is very small. Therefore, the
variation of the pressure jump along the interface is small and not
apparent in the figure. We observe the accurate resolution of the
pressure discontinuity up to the interface without over shoots, the
Gibbs phenomenon.

\begin{table}
\caption{\label{tab:Bubble-Oscillation-2D-1}Bubble Oscillation in $2D$, domain
$[0,2]\times[0,2]$, $\rho=\{1,0.05\}$, $\nu=\{0.0005,0.0005\times0.01\}$,
$\sigma=0.5$, $R_{0}=0.8$, $n=2$ and $\epsilon=0.05$. with theoretical
period $T_{2}=2.56638189$}

\centering{}\begin{tabular}{|c|c|c|}
\hline 
Mesh Size  & Period  & Error in Percentage\tabularnewline
\hline
\hline 
20x20  & 3.4602 & 17.45\tabularnewline
\hline 
40x40  & 3.0476 & 9.95\tabularnewline
\hline 
80x80  & 2.8962 & 6.38\tabularnewline
\hline 
160x160  & 2.7643 & 3.93\tabularnewline
\hline
\end{tabular}
\end{table}

\subsubsection{Bubble Oscillation in Three Dimension}

For a $3D$ droplet under zero-gravity, when the initial position
of the droplet interface with small perturbation \cite{SussmanSmereka1997}
is given by \[
r(\theta,t)=R_{0}+\epsilon P_{n}(cos(\theta)),\]
 where $R_{0}$ is the radius of the droplet, $P_{n}$ is the Legendre
polynomial of order $n$, $\epsilon$ is the amplitude of the perturbation,
then the frequency of the droplet oscillation is given by \[
\omega_{n}=\sqrt{\frac{1}{We}\frac{n(n-1)(n+1)(n+2)}{R_{0}^{3}(n+1+n\lambda)}}\]
 where $\lambda=\frac{\rho_{g}}{\rho_{l}}$, and $We=\frac{\rho_{l}LU^{2}}{\sigma}$
is the Weber number. 

\begin{figure}
\caption{\label{fig:Bubble-Oscillation-in-3D-Velocity-Field}Bubble Oscillation
in $3D$, Velocity Field on a Slice through the Center of the Droplet
at $t=1.8$}

\centering{}\includegraphics[scale=0.3]{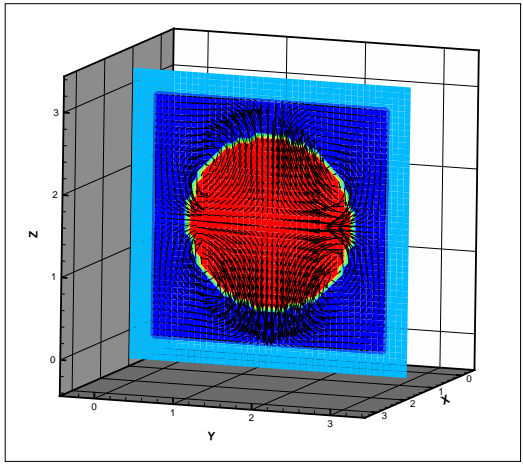} 
\end{figure}

We perform the $3D$ simulation with the domain as $[0,3]\times[0,3]\times[0,3]$.
We use $Re=\frac{\rho_{l}LU}{\mu_{l}}=2000$, $We=\frac{\rho_{l}LU^{2}}{\sigma}=1$
and $\rho_{d}=1$, $\rho_{o}=0.001$, $\mu_{d}=0.0005$, $\mu_{o}=0.000005$.
Table \ref{tab:Bubble-Oscillation-in-3D} shows the convergence of
the oscillation period under mesh refinement. Figure \ref{fig:Bubble-Oscillation-in-3D-Velocity-Field}
shows the velocity field on a slice through the center of the droplet.

\begin{table}
\caption{\label{tab:Bubble-Oscillation-in-3D}Bubble Oscillation in $3D$,
theoretical period is $T_{2}=2.222$}

\centering{}\begin{tabular}{|c|c|c|}
\hline 
Mesh Size  & Period  & Error in Percentage\tabularnewline
\hline
\hline 
20x20x20 & 2.433  & 9.5\%\tabularnewline
\hline 
40x40x40  & 2.300  & 3.5\%\tabularnewline
\hline 
80x80x80  & 2.184  & 1.71\%\tabularnewline
\hline
\end{tabular}
\end{table}

\subsection{EBM for the Elliptic Interface Problem in Cylindrical Coordinates}

We use the method of manufactured solutions to verify our EBM implementation
for the elliptic interface problem in the cylindrical coordinates.
The computational domain is $r\in[1,1.628]$, $\theta\in[0,0.628]$
and $z\in[0.628]$. The interface position is a sphere relative to
the cylindrical coordinates, given as\[
\sqrt{(r-1.314)^{2}+(\theta-0.314)^{2}+(z-0.314)^{2}}=0.2.\]
We use $\rho_{d}=0.811$ for the density inside the sphere and $\rho_{o}=1.03$
for the density out side. We solve the equation\begin{equation}
\nabla\cdot\frac{\nabla p}{\rho}=f\label{eq:elliptic-equation}\end{equation}
where $f$ is a given function. We use \[
p(r,\theta,z)=e^{-\frac{r^{2}+\theta^{2}+z^{2}}{5}}\]
as the exact solution and substitute into the elliptic equation (\ref{eq:elliptic-equation})
to calculate the right hand side $f$. Table \ref{tab:Mesh-Convergence-Study-Elliptic-Cylindrical}
shows the mesh convergence for the this problem. From the table, we
can see that the method is second order accurate. 

\begin{table}
\caption{\label{tab:Mesh-Convergence-Study-Elliptic-Cylindrical}Mesh Convergence
Study for the Elliptic Interface Problem in Cylindrical Coordinates}

\begin{centering}
\begin{tabular}{|c|c|c|c|}
\hline 
Mesh Size & $L_{\infty}$ Error & $L_{2}$ Error & $L_{1}$ Error\tabularnewline
\hline
\hline 
10x10x10 & 0.00018217 & 0.00002381 & 0.00004517\tabularnewline
\hline 
20x20x20 & 0.00009903 & 0.00000659 & 0.00001244\tabularnewline
\hline 
40x40x40 & 0.00001539 & 0.00000174 & 0.00000327\tabularnewline
\hline 
80x80x80 & 0.00000415 & 0.00000045 & 0.00000084\tabularnewline
\hline
\end{tabular}
\par\end{centering}

\end{table}

\subsection{High Speed Two Phase Couette Flow}

\begin{figure}
\caption{\label{fig:Interface-for-the-Rotation-in-Cylinder}Interface for Two
Phase Couette Flow}

\begin{centering}
\includegraphics[scale=0.3]{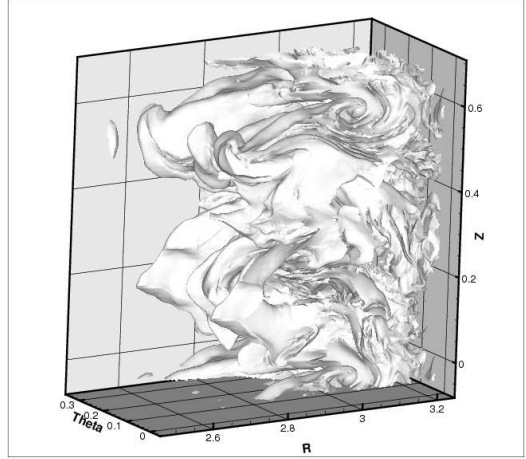}
\par\end{centering}

\end{figure}

Here we use the EBM to simulate high speed two-phase Couette mixing
in a $3D$ angular sector. A more detailed description of this simulation
using the immersed boundary method can be found in \cite{ZhoRayLim12}.
In the current simulation, we use $\rho_{aqu}=1.03g/cm^{3}$ and $\rho_{org}=0.811g/cm^{3}$
for the fluid densities for the aqueous and organic phase respectively,
$\mu_{aqu}=0.0102g/cm\cdot s$ and $\mu_{org}=0.016g/cm\cdot s$ for
the viscosities, and $\sigma=10dyn/cm$ for the surface tension. The
computational domain is an angular sector of the cylinder with $r\in[2.538,3.166]$,
$\theta\in[0,0.314]$, $z\in[0,0.628]$. Figure \ref{fig:Interface-for-the-Rotation-in-Cylinder}
shows the interface position at $t=28\mu s$.

The explanation of the result will be given in another paper. This
example is given here to show the capability of the implemented EBM
to deal with problem with complex interface and of engineering interest.
Under the assumption that there is at most one interface crossing
on each cell edge, there are $2^{8}=256$ possible cases for a mesh
cell. The numerical simulation results for this example show that
all $256$ cases appeared at the same time.

\section{Conclusions}

In this paper, we extended the embedded boundary method to solve the
two-phase incompressible flow. We verify our method by solving the
droplet oscillation problems in $2D$ and $3D$. We also show that
the EBM can be easily extended to solve the elliptic interface problem
in other coordinate systems such as the cylindrical coordinates. Finally,
we simulated the interface contact problem in a rotating cylinder
to show the robustness of the method. Currently, we did not consider
the jump condition due to the discontinuous viscosity. This will be
addressed in the future.

\label{sec:conclusion}

\section{Acknowledgment}

This work is supported in part by Army Research Office W911NF0910306.
Computational resources were provided by the Stony Brook Galaxy cluster,
the Stony Brook/BNL New York Blue Gene/L IBM machine.

\bibliographystyle{plain} \bibliographystyle{plain} \bibliographystyle{plain}
\bibliographystyle{plain}
\bibliography{../refs}

\appendix

\section{Approximate Jump Conditions for Decoupling Velocity}

\label{sub:Appendix-Approximate-Jump-Conditions}

This part is modified from \cite{KangFedkiwLiu2000}. Instead of using
the jump condition (\ref{eq:jump_du2}), we can use any equivalent
jump conditions. Due to the incompressibility condition $\nabla\cdot u=0$,
we have $[\mu\nabla\cdot u]=0$, which can also be written as (Zhi,
et al \cite{TanLeLimKhoo2009}) \begin{equation}
[\mu\frac{\partial u}{\partial n}]\cdot n+[\mu\frac{\partial u}{\partial\tau_{1}}]\cdot\tau_{1}+[\mu\frac{\partial u}{\partial\tau_{2}}]\cdot\tau_{2}=0\label{eq:jump_du3}\end{equation}
 Equations (\ref{eq:jump_du1},\ref{eq:jump_du3}) can be written
in matrix notation in 3D as the following: \begin{equation}
\left(\begin{array}{c}
n\\
\tau_{1}\\
\tau_{2}\end{array}\right)[\mu\frac{\partial u}{\partial n}]+\left(\begin{array}{c}
\tau_{1}\\
n\\
0\end{array}\right)[\mu\frac{\partial u}{\partial\tau_{1}}]+\left(\begin{array}{c}
\tau_{2}\\
0\\
n\end{array}\right)[\mu\frac{\partial u}{\partial\tau_{2}}]+\left(\begin{array}{c}
0\\
\tau_{1}\\
\tau_{2}\end{array}\right)f=0\end{equation}
 Multiplying both sides of this equation from the left by \[
\left(\begin{array}{c}
n\\
\tau_{1}\\
\tau_{2}\end{array}\right)^{T},\]
 we have \begin{eqnarray*}
[\mu\frac{\partial u}{\partial n}]+\left(\begin{array}{c}
n\\
\tau_{1}\\
0\end{array}\right)^{T}\left(\begin{array}{c}
\tau_{1}\\
n\\
0\end{array}\right)[\mu\frac{\partial u}{\partial\tau_{1}}]+\left(\begin{array}{c}
n\\
0\\
\tau_{2}\end{array}\right)^{T}\left(\begin{array}{c}
\tau_{2}\\
0\\
n\end{array}\right)[\mu\frac{\partial u}{\partial\tau_{2}}]\\
+\left(\begin{array}{c}
0\\
\tau_{1}\\
\tau_{2}\end{array}\right)^{T}\left(\begin{array}{c}
0\\
\tau_{1}\\
\tau_{2}\end{array}\right)f & = & 0\end{eqnarray*}
 Since the tangential derivative of the velocity is continuous (\cite{LaiLi2001},
\cite{KangFedkiwLiu2000}), we have \begin{eqnarray}
[\mu\frac{\partial u}{\partial n}]+[\mu]\left(\begin{array}{c}
n\\
\tau_{1}\\
0\end{array}\right)^{T}\left(\begin{array}{c}
\tau_{1}\\
n\\
0\end{array}\right)\frac{\partial u}{\partial\tau_{1}}+[\mu]\left(\begin{array}{c}
n\\
0\\
\tau_{2}\end{array}\right)^{T}\left(\begin{array}{c}
\tau_{2}\\
0\\
n\end{array}\right)\frac{\partial u}{\partial\tau_{2}}\\
+\left(\begin{array}{c}
0\\
\tau_{1}\\
\tau_{2}\end{array}\right)^{T}\left(\begin{array}{c}
0\\
\tau_{1}\\
\tau_{2}\end{array}\right)f & = & 0\label{eq:jump_du_system}\end{eqnarray}

Thus if the surface force $f$ has no or very small tangential components
(as in surface tension force) and $[\mu]$ is very small, it is safe
to approximate the jump conditions of the velocity using \begin{equation}
[\mu\frac{\partial u}{\partial n}]=0\end{equation}
 At an interior point away from the interface, we have $[\mu]=0$
and $f=0$, therefore $[\mu\frac{\partial u}{\partial n}]=0$.

\section{Code Structure and Parallelization}

\label{sec:appendix-code-structure-and-parallelization}

In this section, we briefly introduce the data structure, coding for
dealing with partial cells, and the parallelization method used for
the EBM. For simplicity, we use 2D as example.

The EBM uses a structured Cartesian mesh. For a 2D rectangle domain,
we use a two dimensional array to represent the Cartesian mesh. We
need to store extra geometric information for partial cells (partial
area, cell interface length, center, normal) and partial cell edge
(length, center). We also need to store extra states for the partial
cells and partial cell edges. Since the data size needed for partial
cells with internal interface is much larger than the data size needed
for a regular internal cell, we use an extra pointer to a data structure
for the partial cells (called partial cell data structure later on).
The extra data structure is allocated and deallocated dynamically.
A cell type variable is maintained for each cell to distinguish between
internal, partial, external and boundary cell (cell with external
boundary crossing through). For a moving interface, the cell type
could change dynamically. For the cell edge, the same data structure
is used for both whole and partial cell edges since a partial cell
edge needs only a few extra data storage.

Whenever the interface is changed, the partial cells are identified
and the partial cell data structures are set. The partial cell geometric
information is calculated using cell corner component and cell edge
crossing information using the Marching Cubes algorithm \cite{LorCli87,Wang2012}
(with the assumption that there is at most one crossing at each cell
edge). The geometric information for the partial cell edge is also
set accordingly.

The EBM is a finite volume method. A systematic way of enumerating
the cell internal interface and the partial cell edges is used to
calculate the flux of the differential equations.

Due to the structured Cartesian mesh used, it is very easy to parallelize
the code. There are three types of variables: cell center variables
(cell center velocities, pressure), cell face variables (cell face
velocities), and cell corner variables. The parallelization consists
of three steps: 1) pack; 2) send/receive; 3) unpack. The unpack step
is just the reverse step of the pack step. We first allocate a big
array. Then for each cell to be sent, we first pack the cell type,
and then pack all other variables associated with that cell into the
array. Then we send/receive the single array. The unpack process is
similar with the pack process. For each cell in the buffer zone, we
first get the cell type from the received array, and then unpack the
other variables.
\end{document}